\documentclass[preprint]{elsarticle}

\usepackage[numbers]{natbib} 

\usepackage[left=1in,right=1in,top=1in,bottom=1in,includefoot,includehead,headheight=13.6pt]{geometry}

\usepackage{xcolor}
\definecolor{darkolivegreen}{rgb}{0.33, 0.42, 0.18}
\definecolor{celestialblue}{rgb}{0.29, 0.59, 0.82}

\usepackage{amsmath, amsthm, amssymb}
\usepackage{thmtools, thm-restate}
\usepackage{dsfont}
\usepackage{mathtools}
\usepackage{nicefrac}
\usepackage{graphicx}   
\usepackage{subfigure}
\usepackage{color}
\usepackage{setspace}
\usepackage{lineno}
\usepackage{algpseudocode}

\usepackage{booktabs}

\newcommand{\cM}{\ensuremath{\mathcal{M}}}
\newcommand{\cN}{\ensuremath{\mathcal{N}}}

\newcommand{\cP}{\ensuremath{\mathcal{P}}}

\newcommand{\bE}{\ensuremath{\mathbb{E}}}

\newcommand{\bR}{\ensuremath{\mathbb{R}}}

\newcommand{\bomega}{\boldsymbol{\omega}}

\newcommand{\bv}{\ensuremath{\textbf{v}}}
\newcommand{\bff}{\ensuremath{\textbf{f}}}

\newcommand{\bu}{\ensuremath{\textbf{u}}}

\def\[{\left[}
\def\]{\right]}
\def\<{\langle}
\def\>{\rangle}
\def\({\left(}
\def\){\right)}
\def\[{\left [}
\def\]{\right]}
\def\({\left(}
\def\){\right)}

\newcommand{\norm}[1]{\Vert #1 \Vert}

\DeclareMathOperator*{\argmax}{arg\,max}
\DeclareMathOperator*{\argmin}{arg\,min}

\DeclareMathOperator*{\vspan}{span}

\definecolor{revision}{HTML}{16A085}

\newcommand{\utrue}{\mathbf u_{\text{true}}}
\newcommand{\utilde}{\ensuremath{\tilde{\textbf{u}}}}
\newcommand{\RR}{\ensuremath{\textbf{R}}}

\newcommand{\hside}{\ensuremath{\text{HS}}}
\newcommand{\msmooth}{\ensuremath{\cM_{\text{fast}}}}
\newcommand{\mspike}{\ensuremath{\cM_{\text{slow}}}}
\newcommand{\Thetafast}{\ensuremath{\Theta_{\text{fast}}}}
\newcommand{\Thetaslow}{\ensuremath{\Theta_{\text{slow}}}}
\newcommand{\mslow}{\mspike}
\newcommand{\mfast}{\ensuremath{\msmooth}}
\newcommand{\dist}{\text{dist}}

\newcommand{\inserted}[1]{#1}
\newcommand{\new}[1]{#1} 
\newcommand{\revv}[1]{#1} 
\newcommand{\revvv}[1]{#1}

\newtheorem{remark}{Remark}
\begin{document}

\title{Bias and Multiscale Correction Methods for Variational State Estimation}

\author[1]{F. Galarce\corref{cor1}}
\ead{felipe.galarce@pucv.cl}
\affiliation[1]{organization={School of Civil Engineering. Pontificia Universidad Católica de Valparaíso, Chile}}
\cortext[cor1]{Corresponding author}

\author[2]{J. Mura}
\affiliation[2]{organization={Department of Mechanical Engineering, Universidad Técnica Federico Santa María, Santiago, Chile},}

\author[3]{A. Caiazzo}
\affiliation[3]{organization={Weierstrass Institut für Angewandte Analysis  und  Stochastik,  Leibniz-Institut  im  Forschungsverbund  Berlin  e.V. (WIAS)}, 
state={Berlin}, country={Germany}}

\begin{abstract}
\revvv{Data assimilation performance can be significantly impacted by biased noise in observations, altering the signal magnitude and introducing fast oscillations or discontinuities when the system lacks smoothness. To mitigate these issues, this paper employ} variational state estimation using the so-called parametrized-background data-weak method. 
\revvv{This approach} relies on a background manifold parametrized by a set of constraints, enabling \revvv{the} state estimation by solving a minimization problem on a reduced-order background model, subject to constraints imposed by the input measurements.
\revvv{The proposed formulation incorporates a novel bias correction mechanism
and a manifold decomposition that handles rapid oscillations by treating them as slow-decaying modes based on a two-scale splitting of the classical reconstruction algorithm.}
\revvv{The method is validated in different examples, including the assimilation of biased synthetic data, discontinuous signals, and Doppler ultrasound data obtained from experimental measurements.}
\end{abstract}


\maketitle
\section{Introduction}
Deriving physics-based predictions from limited observations is a multidisciplinary research pursuit, resonating in fundamental science and practical engineering applications. This research finds its systematic exploration in fields such as \textit{data assimilation} and \textit{inverse problems}.
\new{In those fields, the performance of data assimilation algorithms worsens when measurements are corrupted by biased noise, or when the system dynamics lack smoothness, such as in the presence of fast oscillations or discontinuities. This paper main methodological contributions gravitates towards this very issue.}

Data assimilation techniques play an indispensable role in various scientific and engineering domains \cite{carrassi2018,jin2018,cheng2023}. For instance, accurately predicting weather conditions requires integrating data from multiple sources, including satellite observations, ground-based weather stations, and computer simulations. Data assimilation techniques can help \new{to} combine these diverse datasets to improve the accuracy of weather forecasts, understanding and predicting ocean circulation patterns \cite{comse2009data}, refining climate models and merging data from satellites \cite{earthLahoz2014}, among other capabilities like to optimizing reservoir operations and minimize risks associated with water management \cite{ZHANG2021126426}. 

\new{
Data assimilation strategies can be divided into two main categories: state estimation and parameter identification problems \cite{hansen_2007, shaoqing_2020}. The first task relates to computing a time-space resolved field or function, while the second refers to computing the parameters of the underlying model, typically encrypted in a partial differential equation (PDE). Two classical methods for both state and parameter estimation are 4D-var \cite{marco_2020, cstefuanescu2015} and Kalman filters \cite{emmanouil2012, ma11112222}, provided sufficient \textsl{identifiability} \cite{albers_2019}, or sensitivity of the measurements with respect to the parameters to be computed.}
\inserted{
In biomedical imaging, parameter identification plays a relevant role in supporting an image-based diagnosis of tissue and cardiovascular diseases. For example, diffusion tensor imaging (DTI) can be used to estimate the diffusion tensor, a parameter that characterizes the diffusion of water molecules in tissue \cite{DTI}. This information can detect abnormalities in tissue structure, which may indicate disease.
Other applications} have been focused on estimating mechanical tissue parameters from magnetic resonance or ultrasound elastography images. \new{Namely, in \cite{pattison_etal_2014}, the research is focused on the estimation of hydraulic permeabilities, provided an underlying poroelastic model for the tissue. 
Further applications concern the estimation of 
tissue mechanical parameters (elastic moduli, coupling densities, viscoelastic properties) based on magnetic resonance elastography. We refer the reader to \cite{hirsch2017magnetic} for a detailed review.}

Similarly, \new{setting up realistic and patient-specific} blood flow simulations require the estimation of several parameters\new{, which can be inferred, e.g.,} from magnetic resonance images. An example can be found in \cite{lombardi2014}, where sequential parameter estimation was used within the context of hemodynamics, or in \cite{garay_pinns}, which 
followed a deep-learning strategy to accomplish a similar chore. In \cite{chabiniok2012}, the cardiac function is assessed by means of data assimilation techniques and a proper model for the heart's mechanical behavior. For a comprehensive review of computational methods for inverse problems in clinical applications we refer also the reader to \cite{nolte2022}).

State estimation refers to the problem of recovering a full solution (a \textit{state}) on a specific domain of interest, which might differ from the region where the observations are available.
State estimation \revv{strategies} have been used in several applications related to biomedical imaging, such as the estimation of pressure drops in the cardiovascular system \cite{SVIHLOVA2016108, bertoglio_pdrop}, the estimation of displacement and pressure fields from brain elastography data \cite{galarce-etal-2023}, or the the estimation of velocity fields from ultrasound Doppler images \cite{GGLM2021}. \new{The later has even been extended for the estimation of medical biomarkers \cite{GLM2021} and an additional strategy has been proposed to deal with inter-patient geometric variability \cite{GLM2021_2}. These results with synthetically generated data calls for a methodological improvement to deal with biased, noisy or partially corrupted data. This is the main topic of this research.}

If the data assimilation procedure involves the numerical solution of an underlying PDE, the computational cost can become prohibitive\new{, in the case when the forward model needs to be
solved multiple times, for} the increasing the size of the underlying discretization, \revvv{or due to the large size
of the state variable to be sought.}
To overcome those issues, a further key ingredient of modern data assimilation \revv{methods} is the usage of reduced-order models (ROMs), i.e., constructing low-dimensional solution spaces tailored to the physical problem of interest \cite{benner2017model, nicola_dl_mesh}. 

\revvv{
This paper focuses on \revvv{variational} data assimilation based on the so-called Parametrized Background Data-Weak (PBDW) method \cite{MPPY2015}.
Variational approaches addresses the data assimilation as a minimization problem for a cost 
functional which accounts for an underlying model (described by a partial differential equation) and for available knowledge on the considered problem (observations, noise, parameters). Depending on the particular formulation, and for the confidence on the model and the available data, variational method can be formulated balancing between these two aspects and including both the knowledge on the phsyics and available data in the considered cost functional, as optimization contraints, or as additional regularization terms.
The utility of variational strategies extends across numerous domains, exemplified by their successful deployment in various contexts \cite{MPPY2015, evensen2009data}. 

The PBDW approach} is based on computing first a manifold of solution varying the model parameters (a \textit{parametrized background}), which is then approximated
with a linear, low-dimensional, subspace, using ROM techniques (such as the POD -- Proper Orthogonal Decomposition).
The solution to the state estimation problem is then obtained by solving a minimization problem on this \revvv{physics-informed} reduced-order space, combined with a correction depending on available measurements. 

\new{The usage of a low dimensional, \textit{physics-informed} space makes this approach computationally efficient and particularly suitable for data assimilation in the context of applications with limited measurements, i.e., where the reconstruction is sought far beyond the acquisition region. 
In fact, the well-posedness of the reconstruction problem only depends on the relation between the reduced-order space that approximates the solution manifold and the
Riesz representers of the linear functional defining the measurements, and additional regularizations are not required (for a detailed review of the numerical properties of the PBDW, including stability and accuracy, we refer the reader, e.g., to \cite{BCDDPW2017}).
The PBDW has been recently applied in the context of nuclear energy \cite{ARGAUD2018354}, heat transfer problems \cite{heat_pbdw}, and medical imaging, handling assimilation problems for \textit{synthetic} 
ultrasound data \cite{GGLM2021,GLM2021_2} and magnetic resonance images \cite{galarce-etal-2023}.} 
Moreover, in \cite{cohen2022_nonlinearSpaces}, an extension for non-linear sub-spaces was studied, in which the authors provide certified recovery bounds for inversion procedures based on nonlinear approximation spaces. 
\revvv{On the other hand, due to the intrinsic definition of the method, the state estimation via PBDW might suffer from the presence of noise. In fact, the original formulation 
\cite{Maday-pbdw-2015} considered \textit{perfect measurements},
imposing those as constraints for the minimization problem for the state estimation, making its application unsuited in specific contexts. 

To reduce the impact of measurement errors, alternative solutions have been recently proposed, considering observations within a Tikhonov regularization terms (see, e.g., \cite{Taddei2017,maday-taddei-2019}, as well as the 
detailed analysis presented in \cite{GMMT2019} and the recent
extension to time-dependent PBDW reconstruction \cite{Haik_et_al_2023}).
This approach is, however, limited to the case of \textit{unbiased} noise.
Up to our knowledge, extension of the PBDW that are able to 
robustly handle biased measurement noise have not been reported yet.}

Particularly relevant and challenging applications come therefore from ultrasound (US) images, which are prone to low signal-to-noise ratio due to speckle noise combined with anechoic regions and intensity inhomogeneity. There, the efficient and accurate determination of smooth tissue and their interfaces has taken the attention of a whole ultrasound imaging community (see e.g., \cite{CUI2023104431, YUAN2023106837}).

The possibility of defining a suitable approximation of the underlying physics via a reduced-order model is also a key ingredient.
While this aspect allows for an efficient online stage, on the other hand, treating \textit{slow decaying} problems, i.e., problems whose dynamics are characterized by a long Kolmogorov $n$-width \cite{benner2017model,GREIF2019216} and require an arbitrarily large amount of basis functions to represent the physics of the problem correctly, which is highly prohibitive.
Scenarios in which this slow-decaying may become relevant include, phenomena with different scales, including the reconstruction of physical solutions 
with discontinuities. In the latter scenario, methods used for step detection are usually able  to reconstruct the constant parts of the sought solution but may lose fidelity between jumps. 

Motivated by these challenges, the contribution of this work is twofold. Firstly, we propose an extension of the classical PBDW to \revvv{reconstruct the underlying physical state also in presence  of} 
biased measurement noise. The proposed method  defines a correction to the classical PBDW solution 
that exploits the assumption of a priori known physics of the noise structure when incorporating the noise into the state estimation.
Secondly, we propose a strategy to handle a particular class of 
slow-decaying dynamics. Namely, we propose a manifold splitting technique
to tackle the case of physical phenomena with localized discontinuities. 

The rest of the article is organized as follows.
Section \ref{sec:data_ass} introduces the formal setting of the considered data assimilation problems and frameworks, while Section \ref{sec:bias} discusses the extensions to account for biased noise
and Section \ref{sec:multiscale} describes the approach
to handle discontinuities. \new{Test problems for validating the proposed methods are included in each section.
A validation using experimental ultrasound measurements is presented in Section \ref{sec:num_example_us},} while Section \ref{sec:conclusions} draws the conclusions.

\section{General setting of data assimilation problems}\label{sec:data_ass}

\subsection{Preliminaries}
Let us begin introducing a Hilbert space $\left(V,\langle \cdot \rangle \right)$, 
referred to as the \textit{ambient} space, \new{endowed with an inner product $\<\cdot, \cdot \>$ and the induced norm $\norm{\cdot}=\sqrt{(\<\cdot,\cdot\>)}$. 
Moreover, denoting with $X$ a closed subspace of $V$, we also introduce the orthogonal projection of $v \in V$ onto $X$ as follows:
$$
\Pi_{X} v = \argmin_{x \in X} \norm{x - v}.
$$

To formulate the data assimilation problem, we seek a function (the \textit{state}) $\utrue \in V$ 
assuming, as input, a set of (noisy) \textit{observations} $\zeta_1,\hdots,\zeta_m$.
In the above definition, and in the upcoming text, bold variables will denote vector fields and non-bold variables will denote scalar fields.

The observation $\zeta_i$ ($i=1,\hdots,m$) is assumed to be generated by applying a suitable functional 
$\hat{\ell}_i: V \to \mathbb R$ to the unknown state, i.e., 
\begin{equation}\label{eq:measurement}
\zeta_i = \hat{\ell}_i(\utrue) + \psi,\quad i=1,\ldots,m,
    \end{equation}
where $\psi$ models the noise associated with the measurements.
}
\revvv{In practice, the choice of $\psi$ shall reflect the
observational noise in the considered application. At this stage, the particular nature of the measurement noise is not relevant to introduce the method. The upcoming numerical examples (Sections 3--5) will consider applications to the case of 
synthetic measurements with biased Gaussian noise, as well as more complex scenarios involving real velocimetry data with unknown noise distribution.
}

\begin{remark}[Parametrized-background as solution of a PDE]
An example of the considered setting is the case of 
an underlying PDE. In this case,
the ambient space $V$ can be defined as a functional space
containing the solution space. Notice that the setting 
can be defined both as the continuous and at the discrete level, i.e., considering $V$ as a discrete (e.g., finite element) space.
\end{remark}

We assume that the functionals $\ell_i$ are linear and independent. 
Using the Riesz representation theorem, 
each measurement can be thus associated with a representer $\bomega_i \in V$ such that:
\begin{equation}\label{eq:l_i_riesz}
\ell_i(\bu) = \< \bomega_i, \bu\> \quad i=1,\ldots,m.
\end{equation}
Hence, the space of measurement can be formalized by the so-called \textsl{observation space} $W_m = \vspan\{\omega_1,\ldots,\omega_m\}$.
Due to the assumption of independence, it holds $\text{dim}(W_m)=m$.

In this setting, the noiseless part of available measurements on a state $\bu \in V$  can be seen as the projection of the unknown state onto $W_m$. Conversely, a particular element $\bomega^* \in W_m$ can be related, through its coefficients in the basis of representers, to $m$ independent noiseless measurements of a particular state in $V$.

The data assimilation problem can be then formulated via an unknown map
\begin{equation}\label{eq:A-data-assim}
\mathcal A: W_m \to V,\;
\mathbf u^* = \mathcal A(\boldsymbol{\omega}^*),
\end{equation}
which associates, to an element in the observation space, a function (state) in the ambient space.


This type of problem is generally ill-posed. Our approach is to 
close the gap of information by adding a set of constraints based on
an underlying physical model. Namely, we assume
some governing dynamics described by a parametric PDE:
$$
\mathcal P(\bu,\theta) = 0,
$$
where $\theta$ belongs to a set of admissible
parameters $\Theta$.

\begin{remark}[Application to ultrasound imaging]
The present work is motivated by applications of data assimilation in hemodynamics, focusing on the state estimation from Doppler ultrasound images \cite{us2015, VFI_Jensen_Nikolov}.
In this case, the described framework considers a governing parametrized PDE for the fluid flow (Navier--Stokes equations) and the Doppler images as an element in the observation space $W_m$.
\end{remark}

\subsection{The parametrized-background data-weak approach}
In this work, we focus on the Parametrized-Background Data-Weak (PBDW)  approach \cite{MPPY2015}. The information concerning the underlying physical model is encrypted in a \textit{physics-informed} manifold
$$
\cM = \{ \bu \in V;~ \cP(\bu,\theta) = 0;~\theta \in \Theta \subset \bR^p \},
$$
composed of a set of solutions to the underlying dynamics
for a finite set of choices of the model parameters.

To formulate a well-posed data assimilation problem, one then considers a reduced-order space $V_n \subset V$, of dimension $n$, that approximated the manifold $\cM$ with reasonable accuracy, i.e., such that
$$
\dist(u,V_n) = \max_{\bu \in V} \norm{\bu - \Pi_{V_n} \bu}
$$
is small. The reduced-order space $V_n$ can be constructed, e.g., by  means of a proper orthogonal decomposition (POD)
using the snapshots defining the manifold $\cM$. \new{The method considered in this paper and the numerical results shown in the following sections are based on this approach. 

\begin{remark}
It is worth emphasizing that the PBDW is mostly a modular framework, and
other model-order reduction techniques can be used at this stage, such as,  the reduced basis method (see, e.g., \cite{RHP2007}, the (Generalized) Empirical Interpolation Method (see, e.g., \cite{MMT2016, BARRAULT2004667}), orthogonal matching pursuits \cite{GGLM2021}, as well as 
approaches based on nonlinear model-order reduction (see, e.g., \cite{cohen2022_nonlinearSpaces}).
\end{remark}
}

Let $\boldsymbol \omega^*$ denote a representer of a particular set of measurements. The PBDW state estimation 
defines \eqref{eq:A-data-assim} via the solution to the minimization problem
\begin{equation}\label{eq:min-pbdw}
	\begin{aligned}
		\bu^* &= \argmin_{\bu \in V} \frac{1}{2}\norm{\bu - \Pi_{V_n} \bu}^2 \\
		\text{s.t. }& \Pi_{W_m} \bu = \boldsymbol {\omega^*}\,.
	\end{aligned}
\end{equation}
\revv{The approach} \eqref{eq:pbdw_eta} was proposed in \cite{MPPY2015} for the case $\eta(\omega) = \omega$ and further analysed in \cite{BCDDPW2017}.
The well-posedness is guaranteed if the condition 
\begin{equation}\label{eq:pbdw-beta}
\beta(V_n,W_m) = \inf_{\bv \in V_n} \frac{\norm{\Pi_{W_m}\bv}}{\norm{\bv}}>0
\end{equation}
is satisfied.
Moreover, introducing the approximation error of the reduced-order
model
\begin{equation}\label{eq:pbdw-epsilon}
\epsilon_n = \dist \left(\bu, V_n\right) = \max_{u \in \cM} \norm{\bu - \Pi_{V_n} \bu},
\end{equation}
the following bound can be derived \cite{MPPY2015}
\begin{equation}\label{eq:pbdw-bound}
	\norm{\bu - \bu^*} \leq \frac{1}{ \beta\(V_n, W_m\)} \epsilon_n(V_n)\,.
\end{equation}

The interpretation of \eqref{eq:pbdw-bound} is straightforward. On the one hand, it states that 
improving the accuracy of the reduced-order model 
(i.e., reducing $\epsilon_n (V_n)$) also the quality of the reconstruction improves.
On the other hand, it shows that 
the error bound depends on the \textit{angle} between the spaces $V_n$ (reduced-order model) and $W_m$ (observations), meaning that the accuracy of the reconstruction depends on how \textit{observable} the reduced model $V_n$ is. 
One must notice that, from the definition \eqref{eq:pbdw-beta} it follows that, for a fixed set of measurements, the larger the dimension of $V_n$ (i.e., for larger $n$), the 
worse the constant $\beta(W_m,V_n)$ will be.

\section{Bias correction via two-step reconstruction}
\label{sec:bias}
Notice that the state reconstruction computed from \eqref{eq:min-pbdw} has an intrinsic error due to the fact that the constrain on the measurements 
does not account for measurement bias. Hence, it is expected that the quality of the reconstruction deteriorates when
$\norm{\Pi_W \utrue -\omega}$ \new{increases}.

To reduce the impact of measurement bias, we considered a modified \revv{strategy} defined by the following optimization problem:
\begin{equation}
    \begin{aligned}
        \bu_{\eta}^* &= \argmin_{\bu \in V} \frac{1}{2}\norm{\bu - \Pi_{V_n} \bu}^2 \\
		\text{s.t. }& \Pi_{W_m} \bu = \eta(\boldsymbol{\omega}^*)\,, 
	\end{aligned}
	\label{eq:pbdw_eta}
\end{equation}
where the \textit{bias corrector} $\eta: W_m \to W_m$ is a suitable function chosen depending on the given biased observations
$\boldsymbol{\omega}^*$. 

The core idea behind \eqref{eq:pbdw_eta} is to use the
reconstruction $\bu_{0}^*$ from the classical PBDW
problem \eqref{eq:min-pbdw}, which is equivalent to 
$\eqref{eq:pbdw_eta}$ for the special 
case $\eta(\omega) = \omega$, to build a suitable 
approximation of the measurement bias 
and, subsequently, bias corrector.

Namely, we introduce a random noise model
\begin{equation}\label{eq:noise-model-R}
\RR: V \to W_m
\end{equation}
which maps a state in $V$ into noisy measurements.
In practice, $\RR$ can be defined via a random distribution depending on a given set of parameters (e.g., a Normal distribution), or
based on real measurements or simulated data.

\revvv{Notice that the noise model contained in the operator $\RR$ is assumed to be general enough to admit a noise structure -- and a bias -- which depends on the state.}
In general, the form of the noise function shall be driven by knowledge on the application and/or system of interest.
The setting \eqref{eq:noise-model-R} is motivated by the possibility of incorporating noise artifacts (such as bias) of the observation process that depends, e.g., on the magnitude of the measured quantity.

Next, we introduce the discrepancy
\begin{equation}\label{eq:expected-diff}
\xi(\bu;W_m,\RR) := \bE \left( \Pi_{W_m} \bu - \RR(\Pi_{W_m} \bu)  \right) = \Pi_{W_m} \bu - \bE\left[ \RR(\bu )\right]\,.
\end{equation}
The bias corrector is then computed  as
\begin{equation}
	\begin{aligned}
		\eta(\boldsymbol \omega^*) &= 
  \Pi_{W_m} \bu_{0}^* + \xi(\bu_{0}^*;W_m,\RR)  \\
	\end{aligned}
	\label{eq:eta_definition}.
\end{equation}


For a given observation space $W_m$, a set of measurements $\boldsymbol \omega^* \in W_m$, and noise model $\RR$, the bias correction method can be summarized as follows:
\begin{enumerate}
	\item Step 1: Solve \eqref{eq:pbdw_eta} with $\eta_0(\boldsymbol \omega^*) = \boldsymbol \omega^*$ to compute a (biased) reconstruction $\bu_0^*$;
	\item Compute 
 observations $\Pi_{W_m} \utilde$, the error $\xi(\bu_0^*;W_m,\RR)$ 
 using the noise model operator $\RR(\cdot)$;
	\item Step 2: Solve \eqref{eq:pbdw_eta} with 
 the bias corrector $\eta$ computed from \eqref{eq:eta_definition}.
\end{enumerate}

The reconstruction method
defined by the above steps 1--3, as it targets the observation bias, will be hereafter called \textsl{bPBDW}.

\begin{remark}[Well-posedness]
We observe that the two-step reconstruction 
is based on the same spaces used for the 
classical PBDW. Hence, the well-posedness of the original formulation (equation \eqref{eq:pbdw-beta}) ensures also the solvability of the bPBDW optimization problem.
\end{remark}

\revvv{
\begin{remark}[Computational cost of the bPBDW method vs. PBDW] The overall computational cost of the bPBDW reconstruction involves an offline and an online stage. In the offline phase, the solution manifold,  the corresponding reduced order model $V_n$, and the observation space  $W_m$ have to be computed. These steps does not change with respect to the standard PBDW approach. 
However, for the bias-correction method, additional offline computational cost is required for the evaluation of the noise model $\RR(\mathbf u)$, which has to be computed from additional simulation data. The complexity of this calculation depends on the particular problem and on the type of data to be simulated (see, e.g., the Example discussed in Section \ref{sec:num_example_us}).
For the online stage, the bPBDW requires the inversion of 
a $n \times n$ matrix, to solve problem \eqref{eq:pbdw_eta}, thus without any additional cost compared to the
standard PBDW.
\end{remark}}

\new{\subsection{Validation: state estimation from a biased signal}
To validate the proposed approach, we consider the following example (referred in the following as Example 1).
}

Let $\Omega = [0,2\pi]$, let $V = \ell^2(\Omega)$ and let 
us assume that the sought state $\utrue \in V$ belongs to the
set
\begin{equation}
	\cM = \left\{ \bu(x) = A \sin\left(\frac{2 \pi}{T}  x \right), x \in \Omega,~ T\in [\pi, 2 \pi]~A \in [1,2]  \right\}\,.
	\label{eq:simple_manifold}
\end{equation}
This example aims at illustrating the rationale behind the bias correction methodology introduced in this section,
regardless of the physical background. Thus, 
the parametrized background 
is not created from a solution of a PDE,
but defined via a simple sinusoidal function.
\new{The implementation of the numerical methods used for this example are available within the data assimilation software MAD \cite{mad}, which has been
developed and used for different applications based on PBDW (see, e.g., \cite{galarce-etal-2023,GLM2021_2}).}

\begin{figure}[!htbp]
	\centering
	\subfigure[Noise-free measurements]{
		\includegraphics[height=4cm]{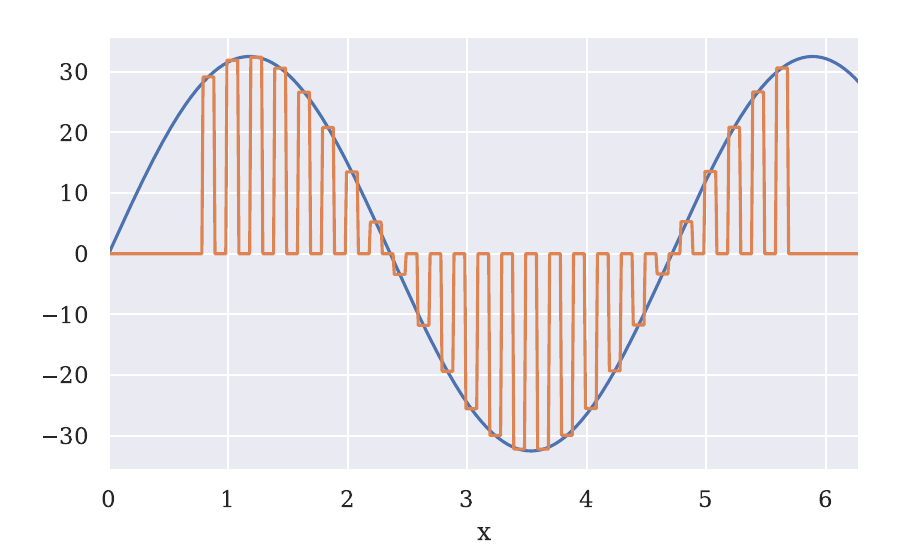}}
	\subfigure[$\omega$ and $\utrue$]{
		\includegraphics[height=4cm]{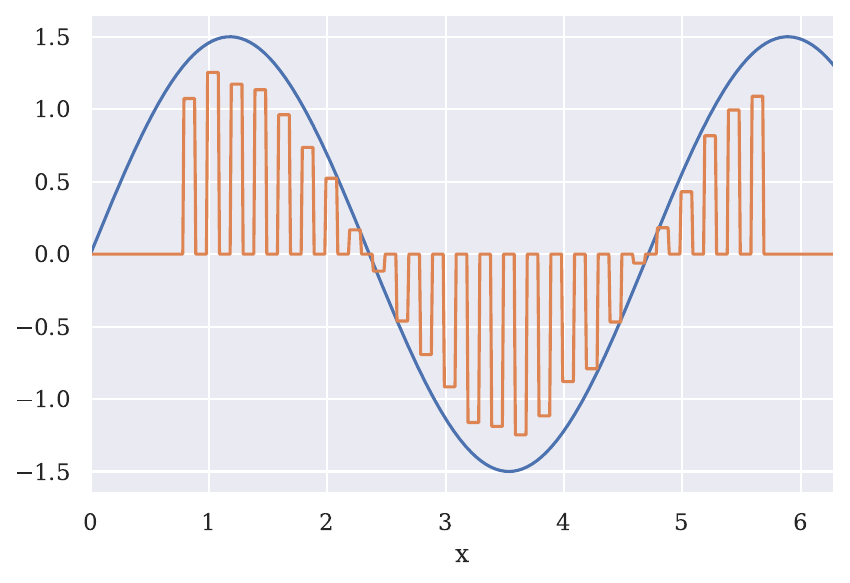}}
	\subfigure[Riesz representers (basis for $W_m$.)]{
		\includegraphics[height=4cm]{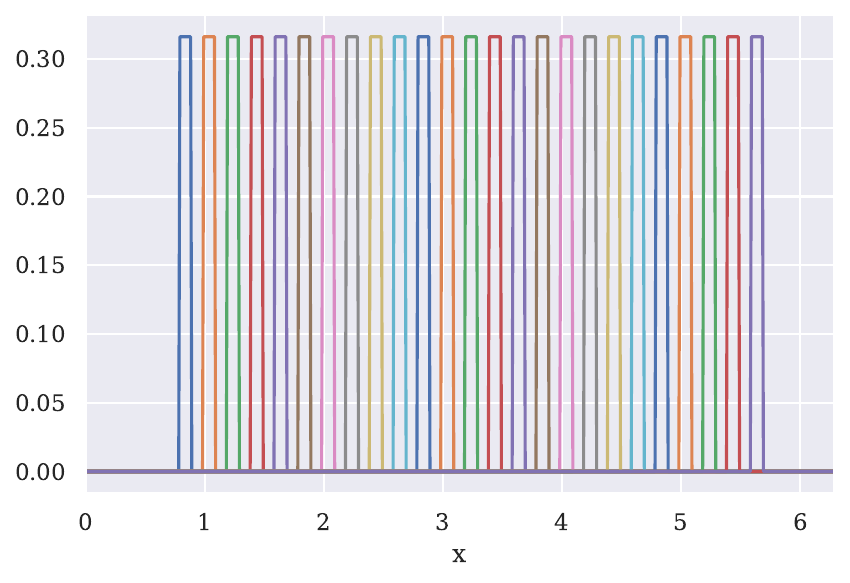}}
	\subfigure[Basis for $V_n$]{
		\includegraphics[height=4cm]{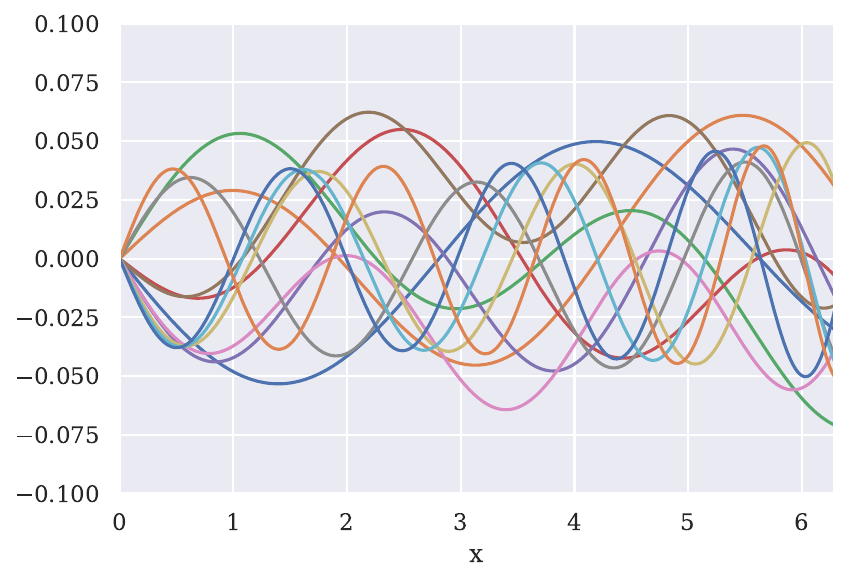}}
	\caption{Setting for Example 1. The plot (a) shows the noise-free measurements, whereas the plot (b) shows the noisy measurements according to the observation model \eqref{eq:test_bias}. Plot (c) and (d) shows the basis for $W_m$ and $V_n$, respectively.}
	\label{fig:u_true_measures}
\end{figure}
Let us  assume that $\cM$ has been obtained sampling uniformly a random distribution in the parameter space, and that the reduced-order model $V_n \subset V$ has been constructed using a POD on the corresponding functions.  

The observation space $W_m$ is built considering 
$m$ equidistant measurement locations.
Moreover, let us consider that the observations are
affected by Gaussian noise $\cN(b(\bu),\sigma)$, where the bias $b(\bu)$ depends on the magnitude of the solution in the form 
\begin{equation}
b(\bu) = \alpha \bu,
\label{eq:test_bias}
\end{equation}
with $\alpha > 0$.  

\begin{figure}[!htbp]
\centering
\includegraphics[height=5cm]{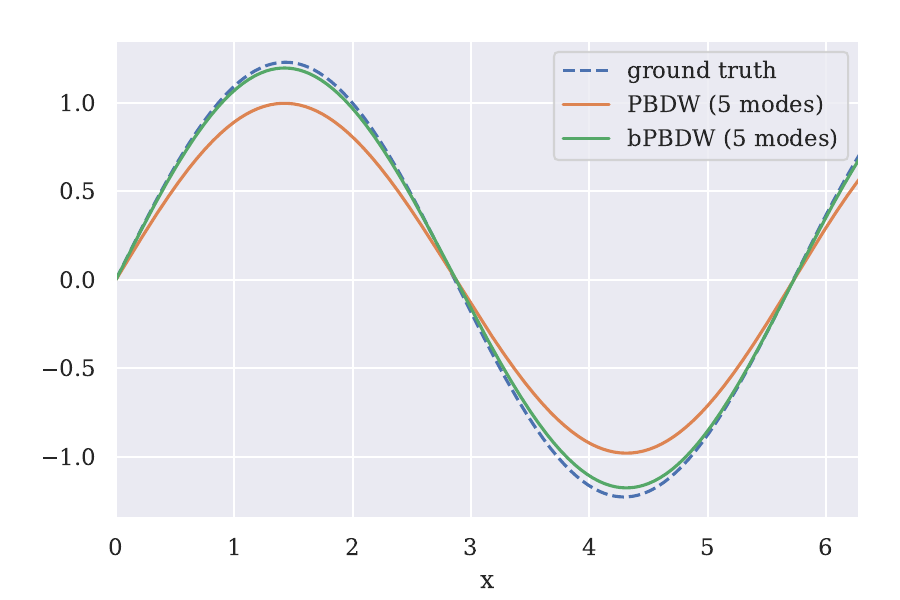}
	\caption{Reconstruction with PBDW and bPBDW for the case of noise parameters $\sigma=A_{gt}/10$  and $\alpha = 0.2$, considering 5 POD modes (optimal dimension for $V_n$, see
 also Figure \ref{fig:simple_experiment_bench}).}
	\label{fig:simple_experiment}
\end{figure}

Figure \ref{fig:u_true_measures} shows a particular
example of a ground truth solution $u_{gt}(x) \in \cM$, for $A = A_{gt} := 32.5$ and $T = T_{gt} := 2\pi$, and the corresponding noisy measurements, for $m=25$. The figure shows the noise-free case to emphasize that the standard PBDW uses these as a model-correction when assimilating the data,
and, therefore,  any strong noise in the measurements might lead to inaccurate results.
In this case, the results for the original reconstruction (Equation \eqref{eq:min-pbdw}) and the bias-corrected one (Equation \eqref{eq:pbdw_eta},
using $\eta(\omega) = \omega$) are
shown in Figure \ref{fig:simple_experiment}, considering 5 POD modes.

\begin{figure}[!htbp]
\centering
\includegraphics[height=4.5cm]{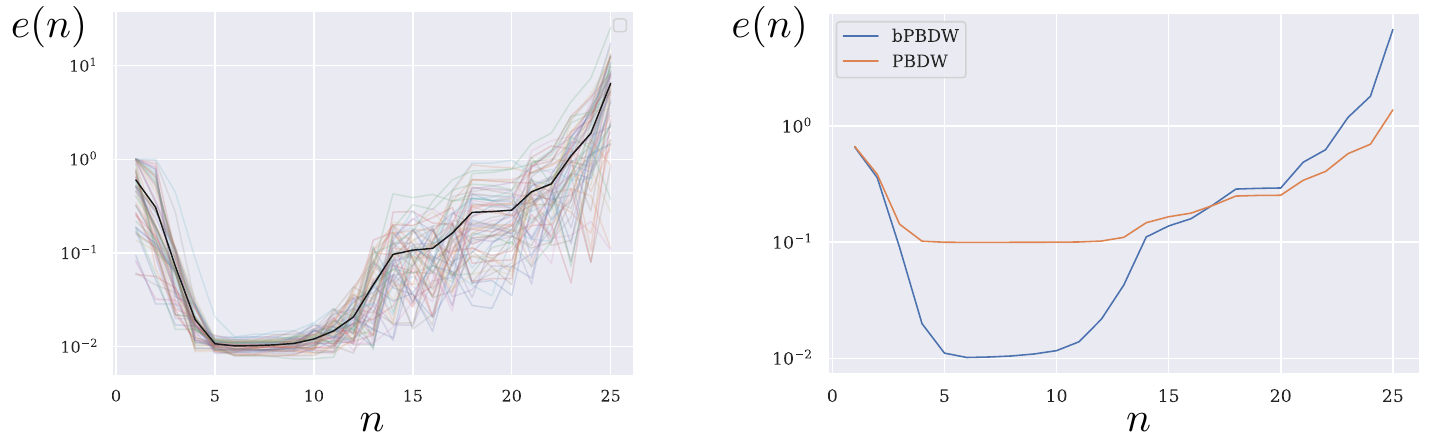}
 \caption{Left: $\ell^2$ relative error $e(n)$ for $\alpha = 0.1$ and $\sigma = A_{gt}/100$ The 64 lighter curves correspond to the error \eqref{eq:error} for every reconstruction, whereas the thicker black curve shows the mean error among the 64 test cases.
 Right: Average reconstruction errors
 for the standard PBDW and for the bias-correction method bPBDW.
 }
 \label{fig:simple_experiment_bench}
\end{figure}

\revvv{For a given ground truth state $\bu_{gt}$, and its reconstruction $\bu^*$, we consider the following metric to assess the reconstruction quality a posteriori:
\begin{equation}
e(n) = \frac{\norm{\bu^* - \bu_{gt}}_{\ell^2(\Omega)} }{\norm{\bu_{gt}}_{\ell^2(\Omega)}},
\label{eq:error}
\end{equation}
evaluated for 64 different ground truth solutions, within the parameter space of the solution manifold range, but different from the ones used for building the ROM. 
Figure \ref{fig:simple_experiment_bench} (left)
shows the behavior of the error \eqref{eq:error} as a function of the dimension of the reduced model $n$ for the 64 reconstructions, as well as the average curve. 

These curves follow the expected behavior of the
error bound \eqref{eq:pbdw-bound}.
In fact, increasing the size of the reduced order model
yields a lower approximation error $\epsilon_n$ \eqref{eq:pbdw-epsilon}
but, at the same time, worsen the stability constant
$\beta(V_n,W_m)$ \eqref{eq:pbdw-beta}. 
In this particular case, the optimal choice
for the size of the reduced-order model is
$n=5$. For this choice, the worst reconstruction error, among the considered 64 test cases, is 5\%.

Figure Figure \ref{fig:simple_experiment_bench} (right) shows the reconstruction errors for the classical PBDW and for the bias-correction method. In this example, the bPBDW is able to achieve reduce the size of the error by  one order of magnitude.
}

\begin{figure}[!htbp]
\centering
\includegraphics[height=4.5cm]{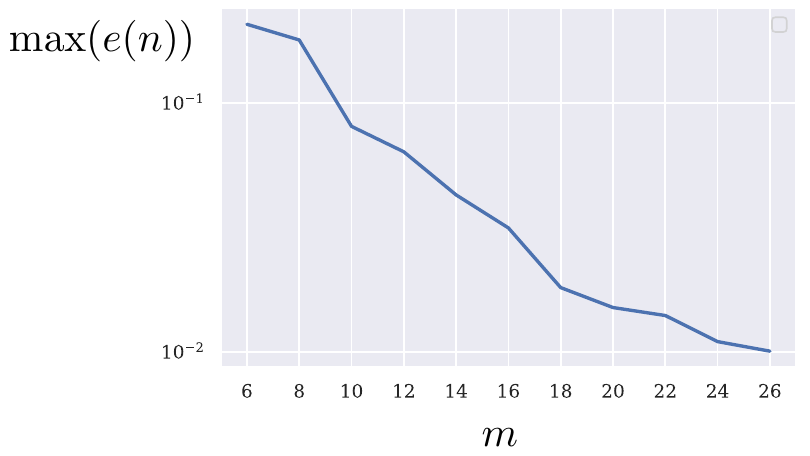}
\caption{Maximal $e(n)$ for bPBDW with respect to the number of measurements $m$ considered, i.e., the dimension of the space $W_m$ ($\alpha = 0.1$ and $\sigma = A_{gt}/100$).}\label{fig:example-1-m-sens}
\end{figure}

\begin{figure}[!htbp]
    \centering
    \includegraphics[height=4.5cm]{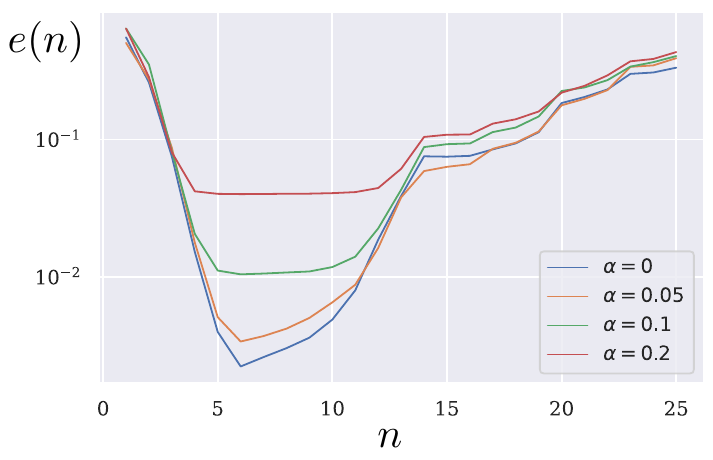}
    \caption{Average $e(n)$ with bPBDW for different bias intensities $\alpha$ (see equation \eqref{eq:test_bias})}
    \label{fig:simple_experiment_bench_comparison}
\end{figure}

\revvv{For the case $n=5$, we performed a sensitivity study with respect to the number of measurements $m$.
Notice that, for a fixed dimension of the reduced space $V_n$, the stability constant $\beta(V_n, W_m)$ increases with $m$
\eqref{eq:pbdw-beta}. The bound on the approximation error
\eqref{eq:pbdw-bound} is therefore expected to improve,
a result which is reflected in Figure \ref{fig:example-1-m-sens},
which shows that the maximal error among the test cases, decreases as the number of sensors (measurement locations) increases.

For the sake of completeness, we observe that, for the noise-free scenario, a perfect reconstruction is obtained for $m \rightarrow \infty$ (setting of unlimited observations setting discussed in \cite{MPPY2015}). 

\begin{remark}[Computational cost associated with additional measurements]
    Adding new measurements might also introduce
    additional computational costs. However, this is
    solely related to the computation of the 
    Riesz representers by solving $m$ problems \eqref{eq:l_i_riesz}.
    In particular, if one considers $m' > m$, measurements, the additional computation of $(m'-m)n$ inner products should be considered when assembling the related normal equations \cite{BCDDPW2017}. Nevertheless, if the sensors are placed at an offline phase and the position does not depend on time, the computational cost related to the the computation of the measurement space $W_m$ only affects the offline stage. The online stage of the data assimilation concerns always  the inversion of a $n\times n$ system of equations. 
\end{remark}
}
\revvv{Finally, Figure \ref{fig:simple_experiment_bench_comparison} compares the reconstruction errors the bias correction method bPBDW for different values of the bias intensities $\alpha$. In the unbiased case ($\alpha=0$) one obtains an optimal reduced-order dimension $n=6$ with the lowest error. As expected, the reconstruction error increases with the bias intensity.
}

\section{Parametrized-background data-weak method for non-smooth dynamics}\label{sec:multiscale}

The quality and the performance of the PBDW reconstruction  \eqref{eq:min-pbdw} depend on the approximation properties of the ROM.
Increasing the dimension of the reduced-order model decreases the inf-sup constant \eqref{eq:pbdw-beta} and, at the same time, increases the overall computational costs of the method.
As a consequence, the approach \eqref{eq:min-pbdw} is not suited in the case of \textit{slow-decaying} dynamics, i.e., 
where high number of POD modes are required to obtain a satisfactory approximation of the physics-informed manifold.

This section introduces a modification
of the original method  \eqref{eq:min-pbdw} to handle
a particular case of slow-decaying dynamics, namely, the presence of 
discontinuities (i.e., very high frequencies) in the
solution, with the goal of obtaining, also in this situation, high-quality reconstructions.

\subsection{Decomposition of the parametrized background}\label{ssec:sPBDW}
\new{Our approach is based on a classical assumption of multiscale splitting, i.e., that the dynamics of the system of interest is characterized by two, ideally well-separated, time and/or space scales.
In the PBDW context, this setting can be formalized assuming that
\begin{itemize}
\item[(i)] we can decompose the manifold $\cM$ as 
\begin{equation}\label{eq:manifold-decomp}
\cM = \mfast + \mslow,
\end{equation}
i.e., that the snapshots can be written as a sum of a smooth, \textit{fast-decaying}, component and a \textit{slow-decaying} components, and
\item[(ii)] it holds
$$
\mfast \cap \mslow = \emptyset\,,
$$
i.e., the scales are fully separated.
\end{itemize}
}

Our approach relies on the fact that, for $\mfast$, we can 
define a suitable reduced-order model $V_n$ to 
approximate the dynamics and we can use \eqref{eq:pbdw_eta} to
obtain an estimation of the fast decaying state component, while
the discontinuities (i.e., the slow decaying modes)
are in $\mslow$.
Notice that, to apply the state estimation \revv{method} considering
different components of the dynamics, it is necessary to handle a decomposition similar to \eqref{eq:manifold-decomp} also for the
available measurements and for the sought state $\bu_{\text{true}}$.
Hence, to deal with the slow decaying modes in $\bu_{\text{true}}$, we seek the closest element to $\bu_{\text{true}}$ in $\mslow$, with the purpose of \textit{capturing} where the discontinuities are located. 

In practice, let us assume to be given set of measurements $\boldsymbol{\omega}^* \in W_m$.
The \revv{proposed strategy} for the reconstruction consists in the following steps.

\textbf{Step 1. Decomposition of measurements}: 
For simplicity, let us assume that the dynamics has a single discontinuity. 
First, we seek a function $\bff^* \in V$ (a \textit{smoother}) using an orthogonal search (OS) criterion to reduce the residual
of the measurement:
\begin{equation}\label{eq:slow-decay-f}
\bff^* = \left( \argmin_{\alpha \in \bR }{\norm{\boldsymbol{\omega}^* - \alpha \Pi_{W_m} \bu_{\text{OS} }}}\right)\frac{ \bu_{\text{OS}}}{\norm{\bu_{\text{OS}}}},
\end{equation}
where
\begin{equation}\label{eq:slow-decay-u-OS}
	\bu_{\text{OS}} = \argmax_{\bv \in \mspike}{ 
 \left\langle \boldsymbol{\omega}^*, \frac{\Pi_{W_m} \bv}{\norm{\Pi_{W_m} \bv}}\right\rangle}.
\end{equation}

\textbf{Step 2. Reconstruction of fast-decaying dynamics}: The smoother is used to apply the PBDW \eqref{eq:pbdw_eta} restricted to the
reduced-order space $V_n$, i.e., considering the \textit{smoothed} measurements 
\begin{equation}\label{eq:slow-decay-omega-f}
\boldsymbol{\omega}_f = \boldsymbol{\omega}^* - \Pi_{W_m} \bff^*
\end{equation}
and computing a corresponding fast-decaying state $\utilde_f \in V_n$.

\textbf{Step 3. Multiscale reconstruction}: To obtain a state reconstruction that accounts for the
slow decaying component in the dynamics, as well as for the presence of measurement bias, we introduce a bias-corrected smoother 
$\bff_u$:
\begin{equation}\label{eq:slow-decay-f-u}
\bff_u = \left(\argmin_{\alpha \in \bR} \norm{ \eta(\omega) - \alpha \Pi_{W_m}  \bu_{\text{OS}}} \right)\frac{\bu_{\text{OS}}} {\norm{\bu_{\text{OS}}}}
\end{equation}
defining the multiscale reconstruction as 
\begin{equation}\label{eq:slow-decay-u-star}
    \bu^* = \utilde_f + \bff_u\,.
\end{equation}

In the upcoming numerical examples, we will refer to the \revv{method} defined by the Steps 1--3 above as \textbf{sPBDW}, as it targets slow-decaying dynamics.

The \revv{procedure} (in particular, Step 2), has been described in the case of a single discontinuity. In the more general case of an arbitrary number of discontinuities, 
this approach can be generalized searching for multiples smoother, repeating
iteratively the process in \eqref{eq:slow-decay-f} until the residue is below
a certain threshold, and constructing 
a slow-decay space $V^\text{slow} \subset V$ spanned by these smoothers. 

\subsection{A priori error analysis}
The scope of this section is to study the behavior of the 
reconstruction error $\norm{\utrue-\bu^*}$ when computing $\bu^*$ using the \revv{method}
defined by the steps \eqref{eq:slow-decay-omega-f}, \eqref{eq:slow-decay-f-u}, and \eqref{eq:slow-decay-u-star} detailed in the Section \ref{ssec:sPBDW}. 
In particular, the goal is to understand the impact of decoupling the fast and slow decaying components of the solution manifold.

Let us consider the linear subspace $V_n$ (reduced-order model for the fast-decaying 
component $\mfast$) and a subspace $V^{\text{slow}}$ 
spanned by the (orthogonalized) smoother functions.

Next, we aim at finding a lower bound for the \textit{inf-sup} constant $\beta(V^{\text{slow}} \oplus V_n, W_m)$, analogous of the constant introduced in \eqref{eq:pbdw_eta} for the classical PBDW. 
Using the orthogonality of the considered subspaces we obtain:
\begin{equation}\label{eq:slow-beta}
	\begin{aligned}
		\beta(V^{\text{slow}} \oplus V_n, W_m)^2 &= \inf_{\mathbf v \in  V^{\text{slow}} \oplus V_n} \frac{\norm{\Pi_{W_m} \mathbf v }^2}{\norm{\mathbf v}^2} \\
		& = \inf_{v_s \in V^{\text{slow}},\;v_f \in V_n} \frac{ \norm{\Pi_{W_m} v_s  + \Pi_{W_m} v_f}^2 }{\norm{v_s + v_f}^2} \\
		& = \inf_{v_s \in V^{\text{slow}},\;v_f \in V_n} \frac{ \norm{\Pi_{W_m} v_s}^2  +\norm{ \Pi_{W_m} v_f }^2 }{\norm{v_s + v_f}^2} \\
	& = \inf_{v_s \in V^{\text{slow}},\;v_f \in V_n} 
  \left( 
  \frac{ \norm{\Pi_{W_m} v_s}^2 }{\norm{v_s}^2} \norm{v_s}^2 + \frac{\norm{ \Pi_{W_m} v_f }^2}{\norm{v_f}^2 } \norm{v_f}^2
  \right) \frac{1}{\norm{v_s + v_f}^2} \\
		& \geq \inf_{v_s \in V^{\text{slow}},\;v_f \in V_n} \min\left\{ \frac{ \norm{\Pi_{W_m} v_s}^2 }{\norm{v_s}^2}, \frac{\norm{ \Pi_{W_m} v_f }^2}{\norm{v_f}^2 } \right\} \frac{ \norm{v_s}^2 + \norm{v_f}^2}{\norm{v_s + v_f}^2} \\
		& \geq \inf_{v_s \in V^{\text{slow}},\;v_f \in V_n} \min\left\{ \frac{ \norm{\Pi_{W_m} v_s}^2 }{\norm{v_s}^2}, \frac{\norm{ \Pi_{W_m} v_f }^2}{\norm{v_f}^2 } \right\}  \\
		& =  \min \{ \beta^F, \beta^S \} ^2. \\
	\end{aligned}
\end{equation}
Hence, the lower bound for the inf-sup constant depends on the \textit{lowest} angle between the measurement space and either $V_n$ or $V^{\text{slow}}$. 
	
Using the decomposition of the manifold and the orthogonality of the fast- and slow-decaying subspaces, one can also conclude that:
\begin{equation}\label{eq:slow-epsilon}
	\begin{aligned}
		\epsilon_n:= \dist(\cM,V^{\text{slow}} \oplus V_n)  &=
		 \sup_{\bu\in \cM^S + \cM^F} \norm{ \bu - \Pi_{V^{\text{slow}}} \bu - \Pi_{V_n} \bu } \\
		&= \sup_{\bu = \bu_s + \bu_f \in \cM^S + \cM^F} \norm {\bu_s + \bu_f - \Pi_{V^{\text{slow}}} \left( \bu_s + \bu_f \right) - \Pi_{V_n} \left( \bu_s + \bu_f \right) } \\
		&\leq \sup_{\bu = \bu_s + \bu_f \in \cM^S + \cM^F} \norm {\bu_s - \Pi_{V^{\text{slow}}}  \bu_s} +  \norm{\bu_f - \Pi_{V_n}  \bu_f  } \\
		&\leq \epsilon^{\text{F}} + \epsilon^{\text{S}},
	\end{aligned}
\end{equation}
where we have introduced the approximation errors 
$\epsilon^{\text{F}} = \dist \left( \mfast, V_n\right)$ and $\epsilon^{\text{S}} = \dist\left( \mslow, V^{\text{slow}}\right)$ are introduced.

Thus the following error bound is thus obtained combining \eqref{eq:slow-beta} and \eqref{eq:slow-epsilon}:
\begin{equation}
	\begin{aligned}
		\norm{\bu - \bu^*} \leq  \frac{1}{ \min\{ \beta^F, \beta^S \} } \left( \epsilon^F + \epsilon^S \right).
	\end{aligned}
	\label{eq:bound_multi-scale}
\end{equation}
Equation \eqref{eq:bound_multi-scale} shows how the model reduction qualities accumulate and, in particular, 
that the reconstruction accuracy is dominated by the \textit{worst measurable} sub-space of $\cM$.

\revvv{
\begin{remark}[Computational cost of the sPBDW vs. PBDW].
The reconstruction with sPBDW involves an offline and
an online stages. With respect to the standard PBDW, the offline phase requires sampling of both the fast and slow decaying manifolds. The additional cost of these operations depends on the particular dynamics considered
for the two scales. 
For the online stage, the additional computational cost of the sPBDW, with respect to the PBDW, is solely related
to the orthogonal search in \eqref{eq:slow-decay-f-u}
for the slow-decaying dynamics, whilst the 
reconstruction of the fast decaying component
\eqref{eq:slow-decay-omega-f} has the same complexity as the standard PBDW. 
\end{remark}}

\new{
\subsection{Validation: assimilation of data with discontinuous background}\label{ssec:example2}
}
This section considers the assimilation of signals that
exhibit discontinuities, using sPBDW, described in section \ref{sec:multiscale}. From the perspective of reduced-order modeling, these discontinuities result in slow decaying dynamics (high number of modes required), making the usage of POD unfeasible. 
\new{As for the previous example (Section \ref{ssec:example2}), the implementation of this example is available within the data assimilation software MAD \cite{mad}.}

Let us consider a domain $\Omega = [0,2\pi]$ and a signal in the form:
\begin{equation}\label{eq:example-2-signal}
\bu = \frac{1}{N_f} \sum_{i=1}^{N_f} A_i \sin\left(\frac{2 \pi}{T_i} x + \delta_i \right) + \beta \hside(x'),   
\end{equation}
where 
$$
\hside(x'):=\begin{cases}
1 & x \geq x' , \\ 0 & \text{otherwise}.
\end{cases}
$$ 
stands for the heavy-side function with a jump at $x'$. The following setting will be referred to as Example 2 in the rest of the section.

In \eqref{eq:example-2-signal}, $N_f$ stands for the number of frequencies to be considered, while $A_i$, $T_i$, and $\delta_i$ ($i=1,\ldots,N_f$) for the amplitudes, periods, and phase shifts of each frequency, respectively.  
The manifold is thus defined as 
	$$
	\cM = \{ u(\theta) \in V = \ell^2(\Omega);~ \theta \in \Theta \},
	$$
	with the parameter space $\Theta = \Thetafast \times \Thetaslow$ with
	\begin{equation}\label{eq:theta-fast-slow}
 \Thetafast = \left( A_1,\ldots, A_{N_f}, \delta_1 \ldots,\delta_{N_f}, T, x\right) \subset \bR^{2+2N_f},\;
 \Thetaslow = \left(x', \beta\right)\subset \bR^{2}\,.
	\end{equation}

In this example, we assume a prior knowledge on form of the state to be estimated (Equation \eqref{eq:example-2-signal}), but
without access to the specific parameters.
As mentioned in Section \ref{sec:multiscale}, the challenge 
of using the POD to solve \eqref{eq:pbdw_eta} concern the discontinuity that the heavy-side function provides, leading to an estimated state where the jumps cannot be located. 
We decompose the manifold $\cM = \mfast \cup \mslow$ according to \eqref{eq:theta-fast-slow}, i.e., defining 
\begin{equation}
	\msmooth = \{ u \in V;~ u = \frac{1}{N_f} \sum_{i=1}^{N_f} A_i \sin\left(\frac{2 \pi}{T} x + \delta_i \right);~ \{ A_1,\ldots, A_{N_f}, \delta_1 \ldots,\delta_{N_f}, T, x\}  \in \Thetafast\},    
\end{equation}
(for the smooth part of the solution), and 
\begin{equation}
	\mspike = \{ u \in V;~ u = \beta \hside(x'); ~\{ \beta, x' \} \in \Thetaslow  \}\,.
\end{equation}
	
\begin{figure}[!htbp]
\centering
\subfigure[$\cM = \msmooth \cup \mspike$]{
    \includegraphics[height=4cm]{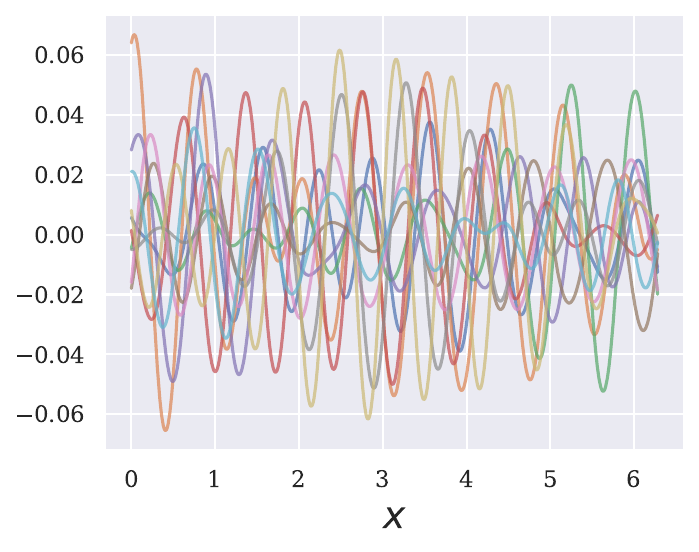}}
\subfigure[$\msmooth$]{
    \includegraphics[height=4cm]{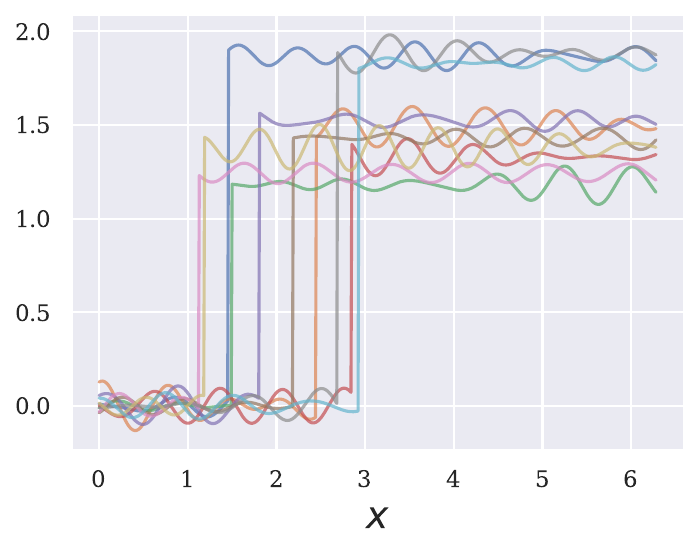}}
\subfigure[$\mspike$]{
    \includegraphics[height=4cm]{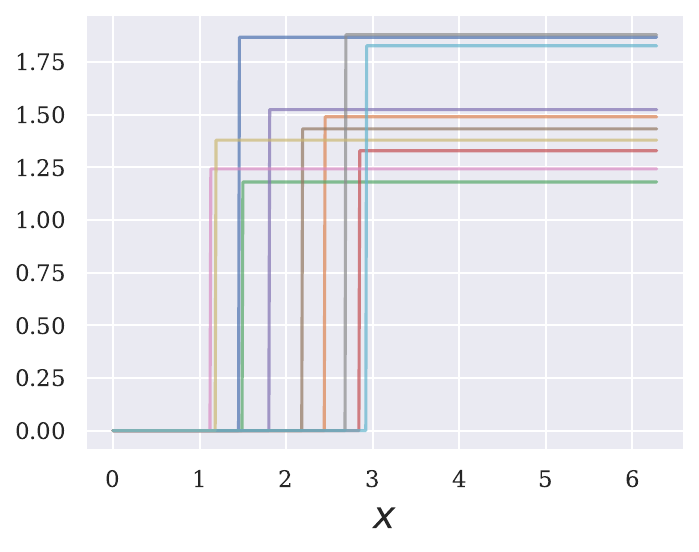}}
\caption{Example 2. Sample of functions of the manifold, considering the complete dynamics (a), the fast decaying component (b) and the discontinuities (c).}
\label{fig:example2-manifolds}
\end{figure}

This splitting allows us to generate a POD reduced-order space $V_n \subset V$ to approximate the smooth part of the snapshots $\msmooth$
(see Figure \ref{fig:example2-manifolds}).

\begin{figure}[!htbp]
\centering
\includegraphics[height=5cm]{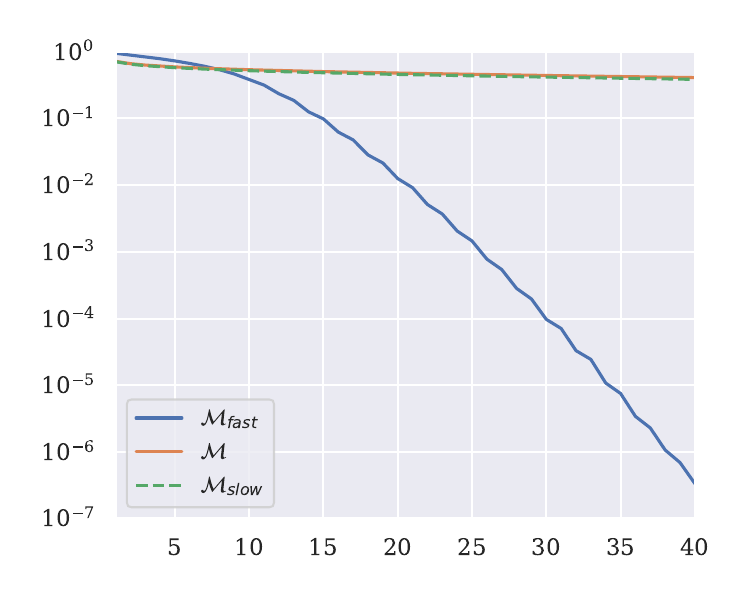}
\caption{Example 2. Approximation error of the reduced-order model
as a function of the reduced dimension $n$.}
\label{fig:example2-error}
\end{figure}

In Figure \ref{fig:example2-error}, we can observe the performance of the model reduction via POD for the different cases. 

In particular, as expected, we observe that  the reduced-order model approximates well the fast-decaying dynamics
$\msmooth$ \revvv{using between 15 and 20 modes}, while the error decay is much slower when considering
the model order reduction of the whole manifold $\cM$. In particular,
Figure \ref{fig:example2-error} shows that, in the latter case,
the approximation error is dominated by the error on the slow
decaying components.

\begin{figure}[!htbp]
    \centering
    \subfigure[$\bu_{\text{true}}$ and $\omega$]{
        \includegraphics[height=5cm]{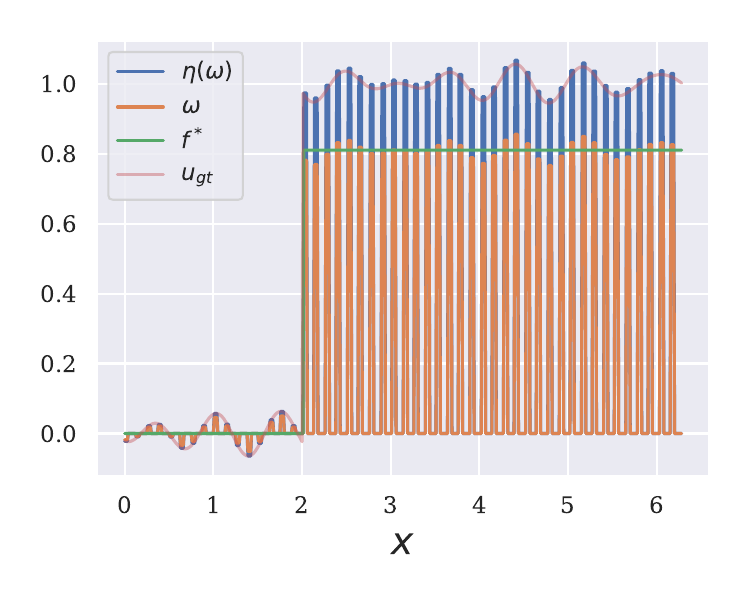}}
    \subfigure[$\bff^*$, $\bff_u$ and $\tilde{\bu}_f$]{
        \includegraphics[height=5cm]{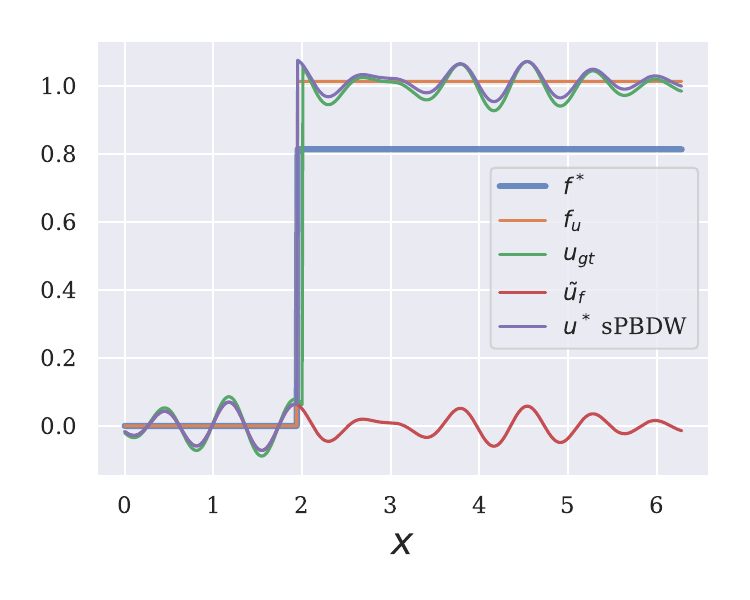}}
    \caption{Example 2.  Figure (a): Example of a ground truth solution, biased measurements ($m=40$), smoother $\bff^*$ and bias corrector $\eta(\omega)$ for the reconstruction with sPBDW. Figure (b): bias-corrected smoother function $\bff_u$ from equation \eqref{eq:slow-decay-f-u}, ground truth solution $\bu_{gt}$ and intermediate solution on fast manifold $\tilde{\bu}_f$. The reconstruction $\bu^* = \tilde{\bu}_f + \bff_u$ is also depicted.}
    \label{fig:step_mani_reconstruction}
\end{figure}

The results of the reconstruction combining the multiscale splitting (sPBDW) and the
bias correction (bPBDW) are shown in Figure \ref{fig:step_mani_reconstruction}. One clearly sees the effect of using the reduced space $V_n$ to reconstruct the smooth dynamic, i.e., once the discontinuity is eliminated via the smoother $\bff_b$ using the orthogonal search. 
	
\begin{figure}[!htbp]
\centering
\includegraphics[height=5cm]{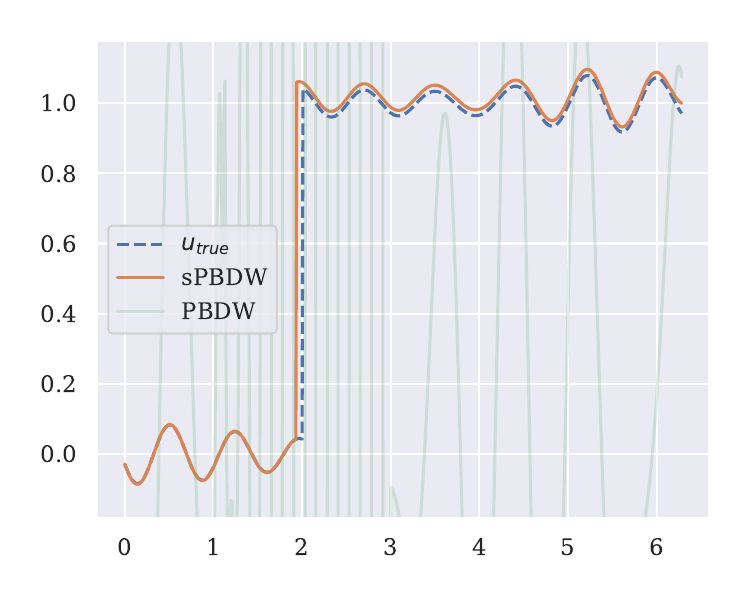}
\caption{Example 2: Comparison of the reconstructions obtained with the classical PBDW (green line) and with 
the multiscale splitting (orange), in the case $n=20$, computed as optimal similarly as it was done for the previous example.}
\label{fig:step_mani_reconstruction2}
\end{figure}

Figure \ref{fig:step_mani_reconstruction2} compares
the signal reconstruction using the classical PBDW and the multiscale approach. We observe that the discontinuities in the state, if included in the creation of the reduced model, leads to a spurious oscillations in the reconstruction. 
The reconstruction obtained after the multiscale splitting overcomes this problem. A similar behavior has been observed independently from the dimension $n$ of the reduced-order model (results omitted).

\section{Assimilation of \revvv{experimental} Doppler ultrasound data}
\label{sec:num_example_us}
The final example is devoted to demonstrate the potential of the bPBDW to assimilate blood velocity data from ultrasound Doppler images \cite{us2015,VFI_Jensen_Nikolov}. \new{Doppler ultrasound is a widely used image modality, based on a piezoelectric transducer that is capable of sending mechanical waves to human tissues and blood, and then capturing a signal that is scattered back by the moving red blood cells inside the blood flow. This allows the computation of a 2D velocity map of the flow \cite{us2013}.}

\begin{figure}[!htbp]
\centering
\subfigure[LUC facilities]{
\includegraphics[height=4.5cm]{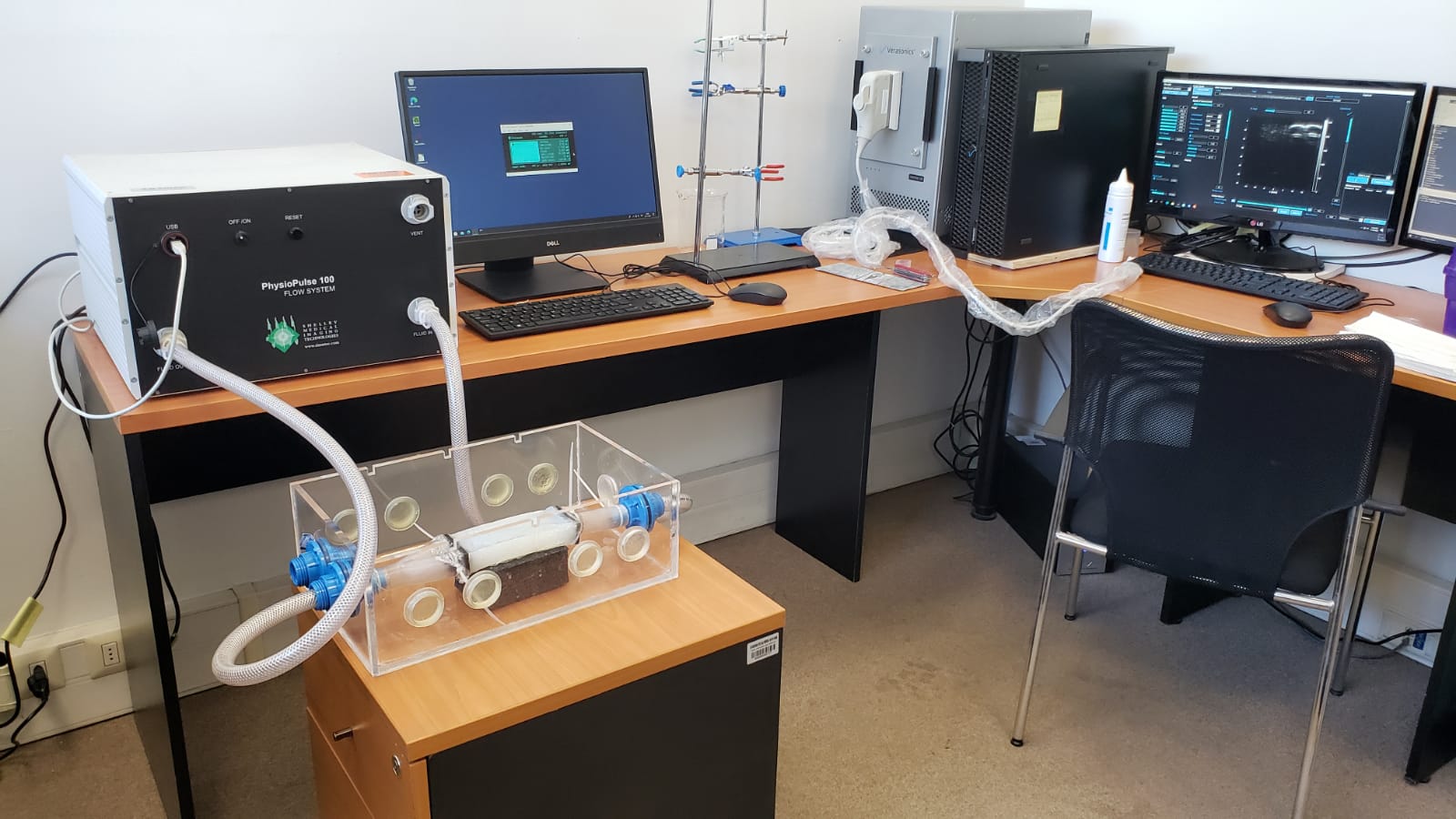}}
\subfigure[Experimental set-up]{
\includegraphics[height=4.5cm]{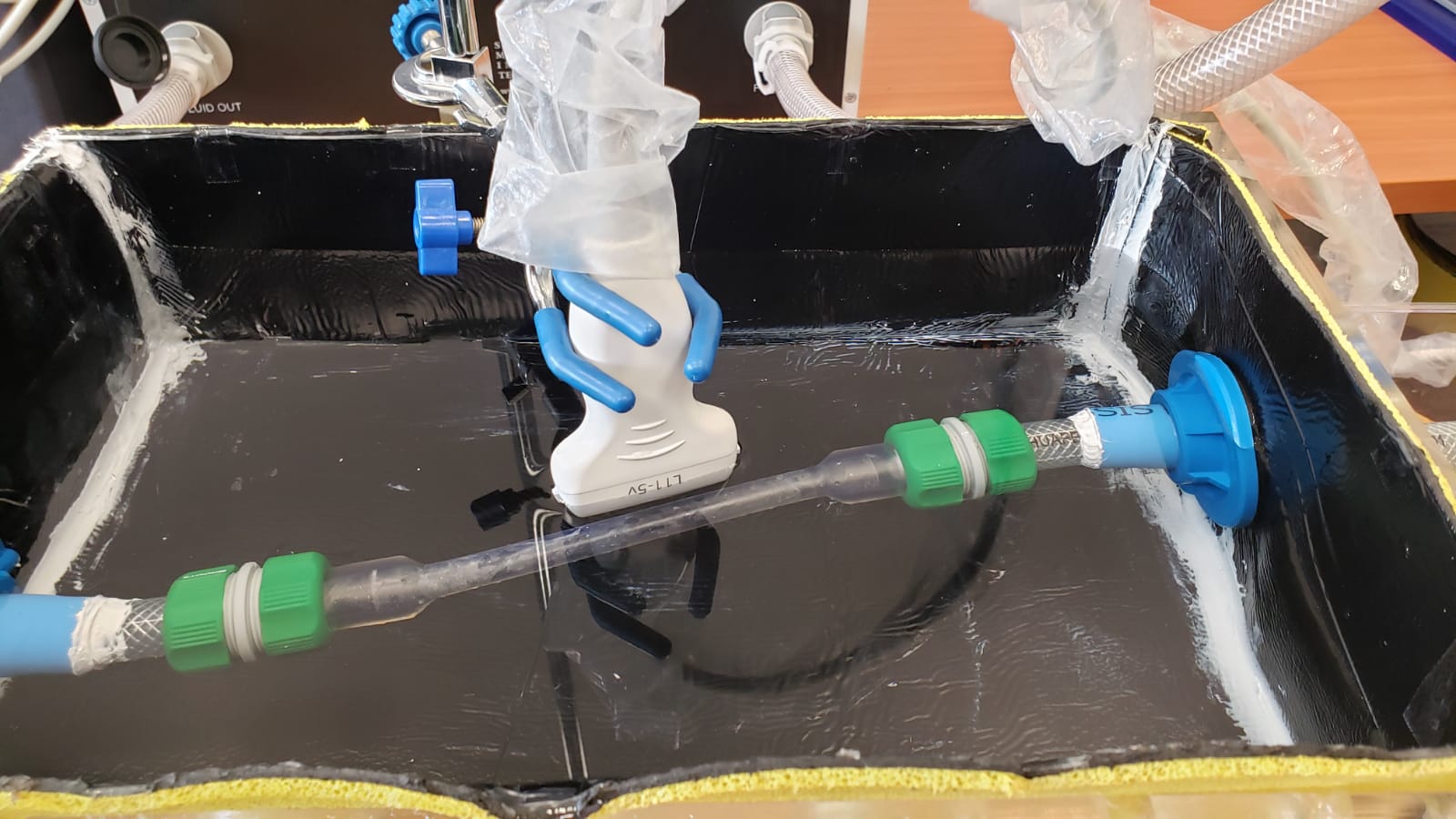}}
\caption{Verasonics 128 machine for Doppler ultrasound measurements at the Laboratorio de Ultrasonido Cuantitativo (LUC), Universidad Técnica Federico Santa María, in Santiago, Chile.}
\label{fig:ex3_luc}
\end{figure}

\begin{figure}[!htbp]
\centering
\subfigure[LUC data]{
\includegraphics[height=6cm]{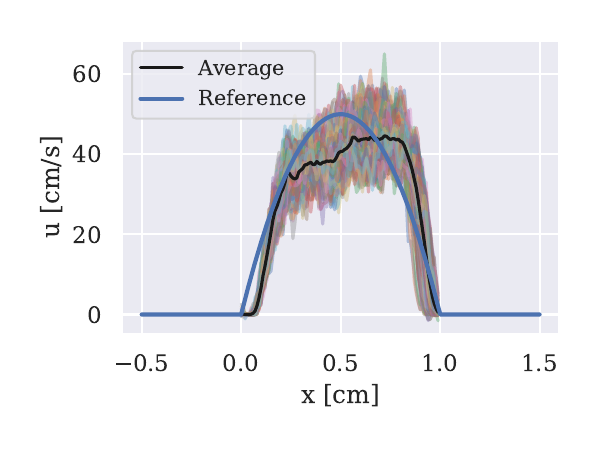}}
\subfigure[Field II data]{
\includegraphics[height=6cm]{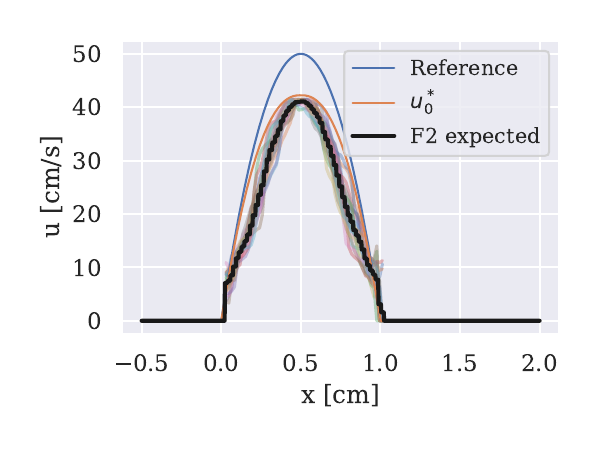}}
\caption{Example 3. (a) Experimental Doppler Flow data acquired with a Verasonics 128 device and and a flow pump Physiopump 100.  b(b) Field II simulated data. (c) Reconstruction with bPBDW.}
\label{fig:ex3_data}
\end{figure}

Namely, our \revv{method} was evaluated by integrating experimental ultrasound data using a color-Doppler technique obtained at a 15-degree insonation angle and a Pulse Repetition Frequency (PRF) of 3.0 kHz. This data was collected at the \emph{LUC} or \emph{Laboratorio de Ultrasonido Cuantitativo USM} (Quantitative Ultrasound Laboratory of the Federico Santa María Technical University). The experimental setup included a tube measuring 20 cm long and 1 cm in diameter, positioned within a regulated setting as depicted in Figure \ref{fig:ex3_luc}. The fluid maintained a steady laminar flow within this setup with a maximum velocity of 50 cm/s. The fluid used was a blood analog suitable for Doppler Ultrasound, characterized by a density of $1037\pm 2\, \text{kg}\,\text{m}^{-3}$ and a viscosity of $4.1\pm 0.1\,\text{mPa}\, \text{s} $. The flow was regulated by the PhysioPulse100 pump, provided by Shelley Medical Imaging Technologies, based in London, Ontario, Canada. A visualization of the acquired data is provided in Figure \ref{fig:ex3_data}(a).

Given the 1-dimensional nature of the phenomena, we focus our attention on the reconstruction of single lines. The flow measurements to be assimilated are obtained with a time-stepping of 7.5 milliseconds.

\begin{figure}[!htbp]
\centering
\subfigure[POD basis for the solution manifold]{
\includegraphics[height=4.5cm]{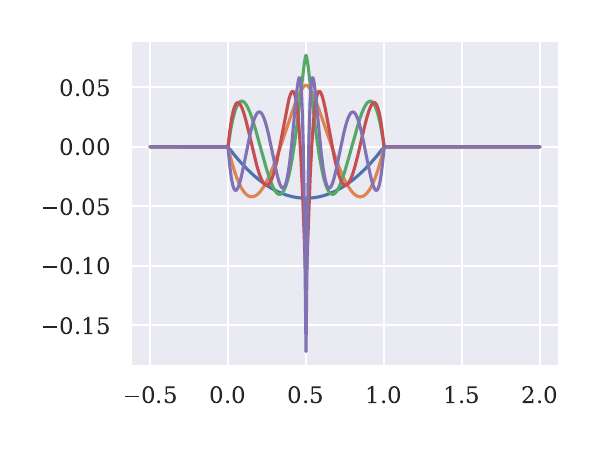}}
\subfigure[$d(u,V_n)$]{
\includegraphics[height=4.5cm]{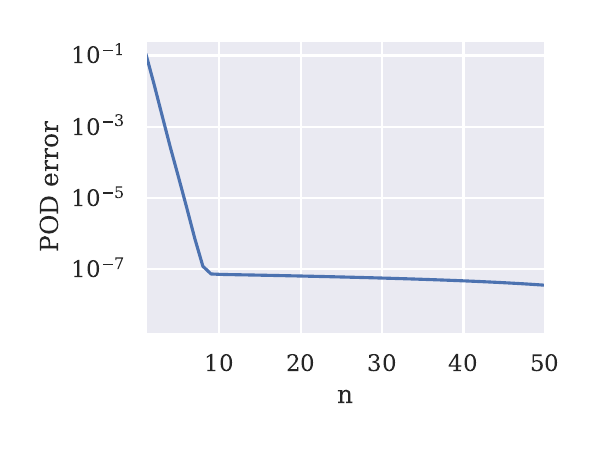}}
\caption{POD of the solution manifold. The Newtonian profile is fully captured with the first parabolic basis, where as the shear-thinning and thickening structures are encrypted in the following modes. The POD error shows a very fast decay, which is a consequence of the simplicity of the dynamics in this particular example.}
\label{fig:ex3_pod}
\end{figure}

The solution manifold for this problem thus contains a collection of functions defined by the set: 
$$
\cM = \left\{ \bu \in V = \ell^2(\Omega) ;~ \bu(r) = v_0 \left[ 1 - \left( \frac{r}{R} \right)^{1+1/\tilde{n}} \right];~ \left( v_0, \tilde{n} \right) \in \Theta \subset \bR^2 \right\}
$$
depending on the parameters $v_0$ (velocity peak), and $\tilde{n}$ (flow behavior index in a power-law rheology). We generate the manifold $\cM$ sampling randomly the parameter space considering 128 values for $v_0$ [cm/s] and $\tilde{n}$, within the ranges $[40,60]$ and $[0.8,1.2]$, respectively. The leading POD modes
for building the reduced basis of $\cM$ and the decay of the POD eigenvalues are shown in Figure \ref{fig:ex3_pod}.
We opt for a non-Newtonian manifold due to experimental and numerical evidence that suggest the relevance of the shear-thinning nature of blood in several similar scenarios \cite{Arzani:2018aa,mella2024_nonNewtonian, mehri_red_2018}. 

Nevertheless, since the experimental fluid used for this particular case is Newtonian, we will hence consider, as target, a 
reference solution $u_{ref}$, with parameters $\tilde{n}=1$, and $v_0=50$ cm/s, i.e., corresponding to the imposed flow boundary condition
up to experimental set-up or measurement errors (blue curve in Figure \ref{fig:ex3_data}).
Note that including the additional viscosity modeling parameter $\tilde{n}$ in the solution manifold will also serve to validate the robustness of the method and its ability 
to reconstruct the appropriate flow pattern, i.e., to identify that the reconstruction corresponds to a parameter $\tilde{n} \approx 1$.

Even if the flow dynamics is  essentially one dimensional, the challenge for the data assimilation lies in the noisy nature of the data.  
The state estimation is attempted both with the standard PBDW\revvv{, which has been successfully employed for the reconstruction from synthetic (not experimental)
Doppler ultrasound data in \cite{GGLM2021,GLM2021},} and with the bPBDW. In both cases, we include the
add-ons proposed in \cite{GLM2021} to handle random noise, consisting in additional restriction is posed on the reconstruction coefficients on the basis for $V_n$, raising from a simple box constraint whose limits are set by the projection of the manifold solutions on the ROM. 

The dimension of the reduced-order space $V_n$ is set  to $n=5$ for all the reconstructions shown in this example.

\begin{figure}[!htbp]
\centering
\subfigure{
\includegraphics[height=6cm]{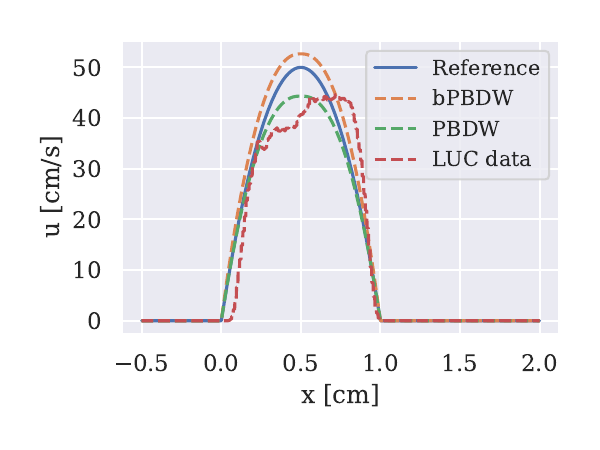}}
\subfigure{
\includegraphics[height=6cm]{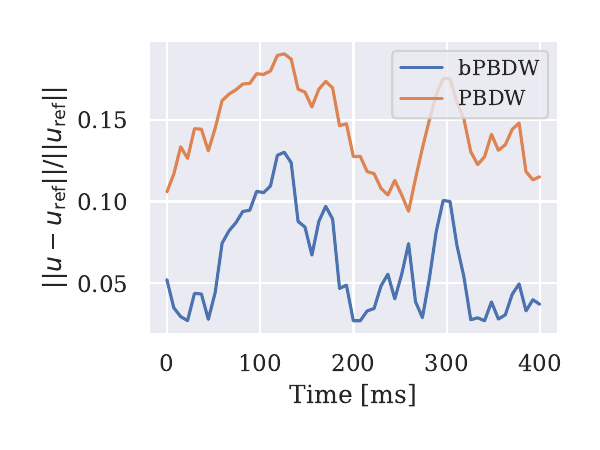}}
\caption{\revvv{Example 3. Left: Reconstruction of average LUC signal comparing the original PBDW and the bias-corrected bPBDW reconstruction with the reference solution $u_{\rm ref}$, for both PBDW and bPBDW. Right: Reconstruction error for the
standard PBDW \cite{GGLM2021,GLM2021} and for the
bias-corrected bPBDW.}}
\label{fig:benchmarck}
\end{figure}

To address the velocity reconstruction, the operator $\RR$ necessary for the correction term for \eqref{eq:eta_definition} 
shall be defined based on an adequate simulation of the ultrasound signal and a preliminary reconstruction $\bu_0^*$. 
For this purpose, we generated synthetic Doppler mappings using the software \textsf{Field II}\footnote{Release 3.30 (April 2021), see: https://field-ii.dk} \cite{field_2_paper,field_2_paper_2}, simulating the signals of the image line by line and then wrapping up contiguous data in time to compute the two-dimensional velocity mapping depicted in Figure \ref{fig:ex3_data}(b).
Namely, \textsf{Field II} has been used to simulate a linear array ultrasound transducer of 196 piezoelectric elements. The simulated device is set to hold similar properties to those in the experimental set-up. 
The simulation produces each line of the ultrasound image and it computes the frequency shift in order to estimate, in each voxel,
the component of the velocity in the direction of the beam. This emulates the scattering of red blood cells groups back to the transducer.

We show the reconstruction over a time-averaged signal on figure \ref{fig:benchmarck}(a). In Figure \ref{fig:benchmarck}(b), the signal in time is reconstructed, showing how bPBDW outperforms PBDW during the whole data acquisition time window.
From a quantitative perspective, testing 40 different input signals, we obtained an average 14.4 \% error using PBDW, and a 6.1 \% error with bPBDW. The error details are depicted in Table \ref{tab:errors}, along with standard deviations, maximal and minimal errors. 
Up to our knowledge, no previous reconstruction of real Doppler data has been made with the original PBDW method, and the improvement
shown in Table \ref{tab:errors} is therefore a major highlight of this work.

\begin{table}[!htbp]
\centering
\caption{Example 3: Overwier of reconstruction errors of PBDW (first row) and bPBDW (second row) using Doppler ultrasound data}
\begin{tabular}{@{}lllll@{}}
\toprule
      & Average(\%)  & Max.(\%)  & Min.(\%)   & Std. dev.(\%) \\ \midrule
PBDW  & 14.4 & 19 & 9.4 & 2.5    \\
bPBDW & 6.11& 13 & 2.7 & 2.7    \\ \bottomrule
\end{tabular}
\label{tab:errors}
\end{table}

\section{Conclusions}\label{sec:conclusions}
\new{This paper investigates and proposes alternative solution methods to overcome two critical issues in handling data assimilation via the Parametrized Background Data Weak (PBDW) method in real scenarios, paving the way for improved accuracy and reliability in estimating physical states from noisy and biased observations. Moreover, the proposed formulations seamlessly integrate into the existing parametrized background framework, allowing for straightforward and efficient implementations. The methods have been demonstrated in simple examples and also validated in the case of experimental Doppler ultrasound data.
}

Firstly, we addressed the challenge of reconstructing the physical state from biased, noisy measurements. The proposed extension to the PBDW (bPBDW) approach employs a two-step procedure to handle the presence of bias. In this procedure, the initial iteration utilizes the classical PBDW solution to approximate the bias correction for the subsequent iteration.
\new{%
We have benchmarked the method assimilating real Doppler ultrasound data using PBDW. Up to our best knowledge, this is an endeavor that has not been reported in the literature yet.
The numerical results show that the proposed approach can qualitatively and quantitatively improve the assimilation of velocity fields from ultrasound data, 
in comparison to classical PBDW.  In particular, the bPBDW bias correction strategy proposed in this work, reduces the average reconstruction error from 14.4\% to 6.1\%, and the maximal error from 19 \% to 13 \%. 
The implications of this enhancement is multilateral, but we can highlight its main relevance on biomedical applications, where blood flow dynamic alterations are inherently related to cardiovascular malfunction  \cite{jing_2022_atherosclerosis, ware_2017_sickle}}.

Next, we investigate the state estimation in the case of discontinuous dynamics. In this case, the discontinuities can be considered as fast-varying scales, introducing 
long Kolmogorov $n$-lengths which are not suitable for the usage of reduced-order models based on POD. 
The proposed approach is based on decomposing the parametrized background in slow (smooth) and fast subsets.
The reduced-order model is used only to recover the smooth component, while an orthogonal search is employed to fit the remaining part.
\new{While the results demonstrate the potential of the proposed manifold splitting, it is worth mentioning that, in selected contexts, this approach could be
linked to different existing two-scale methods, including numerical homogenization \cite{BONIZZONI2024}, or variational multi-scale methods \cite{CASTILLO2019701}, among others. 
Investigating these links will be the object of upcoming research.}

\new{From the methodological perspective, the accuracy gain and the possibility to extrapolate the data beyond the observable region, are the main highlights of the bPBDW strategy. 
On the other hand, it is important to emphasize that the training stage of the reduced model draws a limitation for the method usability in presence of geometric variability.
To overcome the need for patient-specific training, i.e., a large computational cost in the offline phase,  the extension of the proposed approach including additional tools to adapt the pipeline for geometric uncertainties (see, e.g., \cite{GLM2021_2}) represents a promising research direction.
The proposed correction relies on adding an additional reconstruction step to the data assimilation pipeline. While this implies an extra computational cost, several strategies to overcome this issue are also under investigation so that real-time reconstruction can be achieved as done with the original PBDW. 
Further directions also include a data-driven definition of the noise model from clinical data, to reduce 
the overall computational time and target real-time recovery algorithms.}

\section*{Acknowledgments}
We thank Dr. Damiano Lombardi for his valuable insights about model reduction and data assimilation.
F. Galarce sincerely acknowledges the support through research funding VINCI PUCV DI-Iniciación 039.482/2024. J. Mura acknowledge for partial support to the National Agency for Research and Development of Chile (ANID), through the FONDECYT Regular Project \#1230864 and FONDEQUIP mid-scale Scientific and Technological Equipment Funding Program EQM210141. \new{We also thank Eng. Gabriel V\'asquez for his assistance in Doppler acquisition.}

\new{
\section*{CRediT authorship contribution statement}
\begin{itemize}
\item F. Galarce: Conceptualization, Data curation, Validation, Visualization, Software, Formal analysis, Investigation and Writing.
\item J. Mura: Conceptualization, Data curation, Writing and Investigation.
\item A. Caiazzo: Conceptualization, Methodology, Writing, and Formal Analysis.
\end{itemize}
}

\end{document}